\documentclass[11pt]{article}
\usepackage{tikz}
\usepackage[all]{xy}
\usepackage{amsfonts}
\usepackage{amsfonts}
\usepackage{amsfonts}
\usepackage{mathrsfs}
\usepackage{bm}
\usepackage{cite}
\usepackage[colorlinks,linkcolor=blue,citecolor=blue]{hyperref}
\usepackage{amssymb,amsmath} 
\usepackage{makecell}
\usepackage{tikz}
 \allowdisplaybreaks
\usetikzlibrary{arrows,shapes,chains}
 \allowdisplaybreaks

\catcode`!=11
\let\!int\int \def\int{\displaystyle\!int}
\let\!lim\lim \def\lim{\displaystyle\!lim}
\let\!sum\sum \def\sum{\displaystyle\!sum}
\let\!sup\sup \def\sup{\displaystyle\!sup}
\let\!inf\inf \def\inf{\displaystyle\!inf}
\let\!cap\cap \def\cap{\displaystyle\!cap}
\let\!max\max \def\max{\displaystyle\!max}
\let\!min\min \def\min{\displaystyle\!min}
\let\!frac\frac \def\frac{\displaystyle\!frac}
\catcode`!=12

\let\oldsection\section
\renewcommand\section{\setcounter{equation}{0}\oldsection}

\allowdisplaybreaks
\def\pf{\it{Proof.}\rm\quad}

\def\N{\mathbb{N}}\def\Z{\mathbb{Z}}

\def\sn{\sum\limits_{k=1}^n}

\def\su{\sum\limits_{n=1}^\infty}
\def\sk{\sum\limits_{k=1}^\infty}
\def\sj{\sum\limits_{j=1}^\infty}
\def\t{\widetilde{t}}
\def\S{\widetilde{S}}
\def\z{\zeta}

\def\a{^{(A)}}
\def\B{^{(B)}}
\def\C{^{(C)}}

\newtheorem{defn}{Definition}[section]
\newtheorem{thm}{Theorem}[section]
\newtheorem{lem}[thm]{Lemma}
\newtheorem{cor}[thm]{Corollary}

\newtheorem{re}[thm]{Remark}

\begin{document}
\title {\bf Alternating Euler $T$-sums and Euler $\S$-sums}
\author{
{Weiping Wang$^{a,}$\thanks{Email:  wpingwang@yahoo.com, wpingwang@zstu.edu.cn}\quad Ce Xu$^{b,}$\thanks{Email: 19020170155420@stu.xmu.edu.cn}}\\[1mm]
\small a. School of Science, Zhejiang Sci-Tech University,\\
\small  Hangzhou 310018, P.R. China\\
\small b. School of Mathematics and Statistics, Anhui Normal University,\\ \small Wuhu 241000, P.R. China}

\date{}
\maketitle \noindent{\bf Abstract} In this paper, we study the alternating Euler $T$-sums and related sums by using the method of contour integration. We establish the explicit formulas for all linear and quadratic (alternating) Euler $T$-sums and related sums. Some interesting new consequences and illustrative examples are considered.
\\[2mm]


\section{Introduction and Notations}
In our previous paper \cite{XW2019}, we introduced and studied the following two variants of the classical Euler sums,
\begin{align*}
T_{p_1p_2\cdots p_k,q}:=\sum\limits_{n=1}^\infty \frac{h_{n-1}^{(p_1)}h_{n-1}^{(p_2)}\cdots h_{n-1}^{(p_k)}}{(n-1/2)^q}\quad{\rm and}\quad \S_{p_1p_2\cdots p_k,q}:=\sum\limits_{n=1}^\infty \frac{h_n^{(p_1)}h_n^{(p_2)}\cdots h_n^{(p_k)}}{n^q},
\end{align*}
where $q>1,p_1\leq p_2\leq \cdots\leq p_k$ are positive integers, and the quantity $w:={p _1} +  \cdots  + {p _r} + q$ is called the weight and the quantity $r$ is called the degree (or order). We often refer these sums as the the Euler $T$-sums and Euler $\S$-sums, respectively.
Here $h_n^{(p)}$ stands for odd harmonic number of order $p$ defined by
\begin{align*}
h_n^{(p)}:=\sum\limits_{k=1}^n \frac{1}{(k-1/2)^p},\quad h_n\equiv h_n^{(1)} \quad{\rm and}\quad h_0^{(p)}:=0.
\end{align*}
The Euler $T$-sums and Euler $\S$-sums can be seen as variants of classical Euler sums \cite{FS1998}
\begin{align*}
{S_{{p_1p_2\cdots p_k},q}} := \sum\limits_{n = 1}^\infty  {\frac{{H_n^{\left( {{p_1}} \right)}H_n^{\left( {{p_2}} \right)} \cdots H_n^{\left( {{p_k}} \right)}}}
{{{n^q}}}},
\end{align*}
where $H_n^{(p)}$ stands for the $p$-th generalized harmonic number, which is defined by
\begin{align*}
H_n^{(p)}: = \sum\limits_{k=1}^n {\frac{1}{{{k^p}}}},\quad H_n\equiv H_n^{(1)}\quad {\rm and}\quad H_0^{(p)}:=0.
\end{align*}
Like classical Euler sums, the Euler $T$-sums and Euler $\S$-sums can be evaluated by using the method of contour integration developed by Flajolet and Salvy \cite{FS1998}. In \cite{XW2019}, we establish many explicit evaluations of Euler $T$-sums and Euler $\S$-sums via $\log(2)$, multiple zeta and $t$-values. Here, for positive integers $p_1,\ldots,p_k$ with $p_1>1$, the multiple zeta value (MZV for short) \cite{H1992,DZ1994} and multiple $t$-values (MtVs for short) \cite{H2019} are defined by
\begin{align*}
\zeta(p_1,p_2,\ldots,p_k):=\sum\limits_{n_1>\cdots>n_k\geq 1} \frac{1}{n_1^{p_1}n_2^{p_2}\cdots n_k^{p_k}}
\end{align*}
and
\begin{align*}
t(p_1,p_2,\ldots,p_k):&=\sum\limits_{n_1>\cdots>n_k\geq 1\atop n_i \ {\rm odd}} \frac{1}{n_1^{p_1}n_2^{p_2}\cdots n_k^{p_k}}\nonumber\\&=\sum\limits_{n_1>\cdots>n_k\geq 1} \frac{1}{(2n_1-1)^{p_1}(2n_2-1)^{p_2}\cdots (2n_k-1)^{p_k}}.
\end{align*}
As it normalized version,
\begin{align*}
&\widetilde{t}(p_1,p_2,\ldots,p_k):=2^{p_1+p_2+\cdots+p_k}t(p_1,p_2,\ldots,p_k).
\end{align*}
In above definitions of MZVs and MtVs, we put a bar on top of $p_j\ (j=1,\cdots k)$ if there is a sign $(-1)^{n_j}$ appearing in the denominator on the right. Which (one of more the $p_j$ barred) are called the alternating
MZVs, alternating multiple $t$-values. For example,
\begin{align*}
&\zeta(p_1,{\bar p}_2,p_3,{\bar p}_4)=\sum_{n_1>n_2>n_3>n_4>0} \frac{(-1)^{n_2+n_4}}{n^{p_1}_1n^{p_2}_2n^{p_3}_3n^{p_4}_4},\\
&t({\bar p}_1,{\bar p}_2,p_3,p_4)=\sum_{n_1>n_2>n_3>n_4>0} \frac{(-1)^{n_1+n_2}}{(2n_1-1)^{p_1}(2n_2-1)^{p_2}(2n_3-1)^{p_3}(2n_4-1)^{p_4}}.
\end{align*}
In particular, we let
\begin{align*}
{\bar \z}(p):=\sk \frac{(-1)^{k-1}}{k^p}\quad\text{and}\quad {\bar t}(p): = \sum\limits_{k = 1}^\infty  {\frac{{{{\left( { - 1} \right)}^{k - 1}}}}
{{{(k-1/2)^p}}}}\quad (p\geq 1).
\end{align*}
It is clear that the multiple $t$-values can be regard as a level 2 multiple zeta value because of the congruence condition in the summation and of the fact that this value can be written as a linear combination of alternating multiple zeta values.
Recently, Kaneko and Tsumura \cite{KTA2018,KTA2019} also introduced and studied a new kind of multiple zeta values of level two
\begin{align*}
T(p_1,p_2,\ldots,p_k):&=2^k \sum_{m_1>m_2>\cdots>m_k>0\atop m_i\equiv k-i+1\ {\rm mod}\ 2} \frac{1}{m_1^{p_1}m_2^{p_2}\cdots m_k^{p_k}}\nonumber\\
&=2^k\sum\limits_{n_1>n_2>\cdots>n_k>0} \frac{1}{(2n_1-k)^{p_1}(2n_2-k+1)^{p_2}\cdots (2n_k-1)^{p_k}},
\end{align*}
which was called multiple $T$-values (MTVs).

The subject of this paper are alternating Euler $T$-sums and alternating Euler $\S$-sums. First, we give the definitions of alternating harmonic number and odd harmonic number.
Let $p$ and $n$ be positive integers,
the alternating harmonic number ${\bar H}^{(p)}_n$ and odd harmonic number ${\bar h}^{(p)}_n$ are defined by
\begin{align*}
&{ \bar H}_n^{(p)}: = \sum\limits_{k=1}^n {\frac{(-1)^{k-1}}{{{k^p}}}},\quad {\bar H}_0^{(p)}:=0,\quad {\bar H}_n:={\bar H}^{(1)}_n,\\
&{ \bar h}_n^{(p)}: = \sum\limits_{k=1}^n {\frac{(-1)^{k-1}}{{{(k-1/2)^p}}}},\quad {\bar h}_0^{(p)}:=0,\quad {\bar h}_n:={\bar h}^{(1)}_n.
\end{align*}
In the definitions of Euler $T$-sums and Euler $\S$-sums, if replace ``$h^{(p_j)}_n$" by ``${\bar h}^{(p_j)}_n$" in the numerator of the summand, we put a ``bar'' on the top of $p_j$. In particular, we put a bar on top of $q$ if there is a sign $(-1)^{n-1}$ appearing in the
denominator on the right. For example,
\begin{align*}
&T_{p_1{\bar p}_2{\bar p}_3,{\bar q}}=\sum\limits_{n=1}^\infty \frac{h_{n-1}^{(p_1)}{\bar h}_{n-1}^{(p_2)}{\bar h}_{n-1}^{(p_3)}}{(n-1/2)^q}(-1)^{n-1},\quad T_{p_1{\bar p}_2{\bar p}_3p_4,q}=\sum\limits_{n=1}^\infty \frac{h_{n-1}^{(p_1)}{\bar h}_{n-1}^{(p_2)}{\bar h}_{n-1}^{(p_3)}h_{n-1}^{(p_4)}}{(n-1/2)^q},\\
&\S_{{\bar p}_1p_2{\bar p}_3,{\bar q}}=\sum\limits_{n=1}^\infty \frac{h_{n-1}^{(p_1)}{h}_{n-1}^{(p_2)}{\bar h}_{n-1}^{(p_3)}}{n^q}(-1)^{n-1},\quad \S_{p_1{\bar p}_2{p}_3{\bar p}_4,q}=\sum\limits_{n=1}^\infty \frac{h_{n-1}^{(p_1)}{\bar h}_{n-1}^{(p_2)}{h}_{n-1}^{(p_3)}{\bar h}_{n-1}^{(p_4)}}{(n-1/2)^q}.
\end{align*}
The sums of types above (one of more the $p_j$ or $q$ barred) are called the alternating Euler $T$-sums and alternating Euler $\S$-sums, respectively. In \cite{X2020}, we systematic studied all classical (alternating) Euler sums.
In this paper, we study these two above alternating variants of Euler $T$-sums or Euler $\S$-sums by using the methods of contour integration and residue theorem.

Next, we introduce some basic notations. Let $A:=\{a_k\}\ (-\infty < k < \infty)$ be a sequence of complex numbers with ${a_k} = o\left( {{k^\alpha }} \right)\ (\alpha  < 1)$ if $k\rightarrow \pm \infty$. For convenience, we let $A_1$ and $A_2$ to denote the constant sequence $\{1^k\}$ and alternating sequence $\{(-1)^k\}$, respectively.
\begin{defn}\label{def1} With $A$ defined above, we define the parametric digamma function $\Psi \left( { - s;A} \right)$ by
\begin{align}\label{1.1}
\Psi ( { - s;A} ):= \frac{{{a_0}}}{s} + \sum\limits_{k = 1}^\infty  {\left( {\frac{{{a_k}}}{k} - \frac{{{a_k}}}{{k - s}}} \right)}\quad (s\in \mathbb{C}\setminus (\N\cup\{0\})).
\end{align}
\end{defn}
Obviously, if $A=A_1$, then the parametric digamma function $\Psi \left( { - s;A} \right)$ becomes the classical digamma function $\psi \left( { - s} \right)+\gamma$.

\begin{defn}\label{def2} Define the cotangent function with sequence A by
\begin{align}\label{1.2}
 \pi \cot( {\pi s;A} ) &=  - \frac{{{a_0}}}{s} + \Psi ( { - s;A} ) - \Psi ( {s;A})\nonumber \\
  &= \frac{{{a_0}}}{s} - 2s\sum\limits_{k = 1}^\infty  {\frac{{{a_k}}}{{{k^2} - {s^2}}}}.
\end{align}
\end{defn}
It is clear that if letting $A=A_1$ or $A_2$ in (\ref{1.2}), respectively, then it become
\begin{align*}
&\cot( {\pi s;A_1}) = \cot ( {\pi s} )\quad\text{and}\quad \cot( {\pi s;A_2} ) = \csc ( {\pi s} ).
\end{align*}

The Definitions \ref{def1} and \ref{def2} are also introduced in a previous paper \cite{X2020} of the second named author.

\begin{defn}\label{def3} For nonnegative integers $j\geq 1$ and $n$, we define
\begin{align*}
&D\a(j):=\sk \frac{a_k}{k^j},\quad D\a(1):=0,\quad E\a_n(j):=\sn \frac{a_{n-k}}{k^j},\quad E\a_0(j):=0,\\
&{\bar E}\a_n(j):=\sn \frac{a_{k-n-1}}{k^j},\quad {\bar E}\a_0(j):=0,\quad\widehat{E}\a_n(j):=\sn \frac{a_{n-k}}{(k-1/2)^j},\quad \widehat{E}\a_0(j):=0,\\
&\widetilde{E}\a_n(j):=\sn \frac{a_{k-n-1}}{(k-1/2)^j},\quad \widetilde{E}\a_0(j):=0,\\
&\widehat{t}\a(j):= \left\{ {\begin{array}{*{20}{c}} \sk \left(\frac{a_{k-1}}{k-1/2}-\frac{a_k}{k}\right)
   {,\ j=1,}  \\
   {\sk \frac{a_{k-1}}{(k-1/2)^j},\ \ \ \ \ \ \ \;\;\;j>1,}  \\
\end{array} } \right.,\quad {\t}\a(j):= \left\{ {\begin{array}{*{20}{c}} \sk \left(\frac{a_{k}}{k-1/2}-\frac{a_k}{k}\right)
   {,\ j=1,}  \\
   {\sk \frac{a_{k}}{(k-1/2)^j},\ \ \ \ \ \ \ \;\;\;j>1,}  \\
\end{array} } \right.\\
&F\a_n(j)= \left\{ {\begin{array}{*{20}{c}} \sk \frac{a_{k+n}-a_k}{k}
   {,\ j=1,}  \\
   {\sk \frac{a_{k+n}}{k^j},\ \ \ \ \ \;\;\;j>1,}  \\
\end{array} } \right. ,\quad {\bar F}\a_n(j)= \left\{ {\begin{array}{*{20}{c}} \sk \frac{a_{k-n}-a_k}{k}
   {,\ j=1,}  \\
   {\sk \frac{a_{k-n}}{k^j},\ \ \ \ \ \;\;\;j>1,}  \\
\end{array} } \right.,\\
& \widehat{F}\a_n(j)= \left\{ {\begin{array}{*{20}{c}} \sk \left(\frac{a_{k+n}}{k-1/2}-\frac{a_k}{k}\right)
   {,\ j=1,}  \\
   {\sk \frac{a_{k+n}}{(k-1/2)^j},\ \ \ \ \ \ \ \ \ j>1,}  \\
\end{array} } \right.,\quad \widetilde{F}\a_n(j)= \left\{ {\begin{array}{*{20}{c}} \sk \left(\frac{a_{k-n}}{k-1/2}-\frac{a_k}{k}\right)
   {,\ j=1,}  \\
   {\sk \frac{a_{k-n}}{(k-1/2)^j},\ \ \ \ \ \ \ \ \ j>1,}  \\
\end{array} } \right.,\\
&G\a_n(j):=E\a_n(j)-{\bar E}\a_{n-1}(j)-\frac{a_0}{n^j},\quad G\a_0(j):=0,\quad L\a_n(j):=F\a_n(j)+(-1)^j{\bar F}\a_n(j),\\
&M\a_n(j):=E\a_n(j)+(-1)^j {F\a_n}(j),\quad {\bar M}\a_n(j):={\bar F}\a_n(j)-{\bar E}\a_{n-1}(j),\quad n\geq 1,\\
&R\a_n(j):=G\a_n(j)+(-1)^j L\a_n(j),\quad N\a_n(j):=\widehat{E}\a_n(j)+(-1)^j\widehat{F}\a_{n-1}(j),\\
&{\bar N}\a_n(j):=\widetilde{F}\a_n(j)-\widetilde{E}\a_{n-1}(j),\quad S\a_n(j):=N\a_n(j)-{\bar N}\a_n(j)-\frac{a_0}{(n-1/2)^j}.
\end{align*}
\end{defn}

In particular, in Definition \ref{def3}, setting $A=A_1$ or $A_2$ yield
\begin{align*}
&M^{(A_1)}_n(j)=H^{(j)}_n+(-1)^j\z(j),\ {\bar M}^{(A_1)}_n(j)=\z(j)-H^{(j)}_{n-1},\ R^{(A_1)}_n(j)=(1+(-1)^j)\z(j),\\
&M^{(A_2)}_n(j)=(-1)^{n-1}{\bar H}^{(j)}_n+(-1)^j \left\{ {\begin{array}{*{20}{c}} (1-(-1)^n)\log(2)
   {,\ j=1,}  \\
   {(-1)^{n-1}{\bar \z}(j),\ \ \ \ \ \ \;\;\;j>1,}  \\
\end{array} } \right.\\
&{\bar M}^{(A_2)}_n(j)=(-1)^{n}{\bar H}^{(j)}_{n-1}+\left\{ {\begin{array}{*{20}{c}} (1-(-1)^n)\log(2)
   {,\ j=1,}  \\
   {(-1)^{n-1}{\bar \z}(j),\ \ \ \ \ \;\;\;j>1,}  \\
\end{array} } \right.\\
&R^{(A_2)}_n(j)=(-1)^{n-1}(1+(-1)^j){\bar \z}(j),\\
&N^{(A_1)}_n(j)=h^{(j)}_n+(-1)^j\t(j),\quad \t(1):=2\log(2),\\
&N^{(A_2)}_n(j)=(-1)^{n-1}{\bar h}^{(j)}_n+(-1)^j \left\{ {\begin{array}{*{20}{c}} (-1)^n{\bar t}(1)+\log(2)
   {,\ j=1,}  \\
   {(-1)^{n}{\bar t}(j),\ \ \ \ \ \ \ \ \ \quad j>1,}  \\
\end{array} } \right.\\
&{\bar N}^{(A_1)}_n(j)=\t(j)-h^{(j)}_{n-1},\quad \t(1):=2\log(2),\\
&{\bar N}^{(A_2)}_n(j)=(-1)^{n}{\bar h}^{(j)}_{n-1}+\left\{ {\begin{array}{*{20}{c}} (-1)^{n-1}{\bar t}(1)+\log(2)
   {,\ j=1,}  \\
   {(-1)^{n-1}{\bar t}(j),\ \ \ \ \ \ \ \ \ \quad j>1,}  \\
\end{array} } \right.\\
&S^{(A_1)}_n(j)=(1+(-1)^j)\t(j),\quad S^{(A_2)}_n(j)=(-1)^{n-1}(1-(-1)^j){\bar t}(j).
\end{align*}

\section{Lemmas}
In this section, we give some power series expansions for parametric digamma function $\Psi \left( { - s;A} \right)$ and $\cot(\pi s,A)$.

\begin{lem}\label{blem1}(\cite[Theorem 2.2]{X2020}) Let $p$ and $n$ be positive integers, if $|s+n|<1$, then
\begin{align}\label{b1}
\frac{{{\Psi ^{\left( {p - 1} \right)}}\left( { - s;A} \right)}}{{\left( {p - 1} \right)!}}=(-1)^p \sj \binom{j+p-2}{p-1} {\bar M}\a_n(j+p-1)(s+n)^{j-1}.
\end{align}
\end{lem}

\begin{lem}\label{blem2}(\cite[Theorem 2.3]{X2020}) With $\cot(\pi s;A)$ defined above, if $|s-n|<1$ with $s\neq n\ (n\in \Z)$, then
\begin{align}\label{b2}
\pi \cot(\pi s;A)=\frac{a_{|n|}}{s-n}-\sj (-\sigma_n)^j R\a_{|n|}(j)(s-n)^{j-1},
\end{align}
where $\sigma_n$ is defined by the symbol of $n$, namely,
\begin{align*}
\sigma_n:= \left\{ {\begin{array}{*{20}{c}}\ 1
   {,\  n\geq 0}  \\
   {-1,\ n<0.}  \\
\end{array} } \right.
\end{align*}
\end{lem}

\begin{lem}\label{lem2.1} Let $p>0$ and $n$ be a non-negative integer, if $|s-n|<1$ then
\begin{align}
&\frac{\Psi^{(p-1)}(1/2-s;A)}{(p-1)!}=\sj (-1)^{j-1}\binom{j+p-2}{p-1} N\a_n(j+p-1)(s-n)^{j-1},\label{2.1}
\end{align}
and if $|s+n|<1$ then
\begin{align}
&\frac{\Psi^{(p-1)}(1/2-s;A)}{(p-1)!}=\sj (-1)^{j-1}\binom{j+p-2}{p-1} {\bar N}\a_{n+1}(j+p-1)(s+n)^{j-1}.\label{2.2}
\end{align}
\end{lem}
\pf The proofs of this lemma follows the definition of function $\Psi \left( { - s;A} \right)$. \hfill$\square$

Letting $n=0$ in (\ref{2.1}) and (\ref{2.2}) gives
\begin{align}
\frac{\Psi^{(p-1)}(1/2-s;A)}{(p-1)!}=(-1)^p \sj \binom{j+p-2}{p-1} \widehat{t}\a(j+p-1)s^{j-1} \quad (|s|<1).
\end{align}

\begin{lem}\label{lem2.2} Let $m>0$ and $n>1$ be non-negative integer, if $|s-n+1/2|<1$ then
\begin{align}\label{2.4}
\frac{d^m}{ds^m} (\pi \cot(\pi s;A))=(-1)^mm!\sj (-1)^{j-1} \binom{j+m-1}{m} S\a_n(j+m)(s-n+1/2)^{j-1}.
\end{align}
\end{lem}
\pf Lemma \ref{lem2.2} follows immediately from Definition \ref{def2} and Lemma \ref{lem2.1}.\hfill$\square$

From Lemma \ref{lem2.2}, we have
\begin{align}
&\lim_{s\rightarrow 1/2} \frac{d^m}{ds^m} (\pi \cot(\pi s;A))=m!((-1)^m\widehat{t}\a(m+1)-\t\a(m+1)),\\
&\pi \cot(\pi s;A)=\frac{a_0}{s}-2\sj D\a(2j)s^{2j-1}\quad (|s|<1).
\end{align}

\begin{lem}(\cite{X2020})\label{lem2.3} Let $p$ and $n$ be positive integers, if $|s-n+1/2|<1$ and $s\neq n-1/2$, then
\begin{align}\label{2.7}
&\frac{{{\Psi ^{( {p - 1})}}( { 1/2- s;A} )}}{{( {p - 1} )!}}\nonumber\\&=\frac{1}{{( {s - n+1/2} )}^p}\left\{a_{n-1}-\sj (-1)^j\binom{j+p-2}{p-1} M\a_{n-1}(j+p-1)(s-n+1/2)^{j+p-1} \right\}.
\end{align}
\end{lem}
\pf This lemma can be immediately obtained from \cite[Theorem 2.1]{X2020}. \hfill$\square$

If $n=1$ then
\begin{align}\label{2.8}
&\frac{{{\Psi ^{( {p - 1} )}}( { 1/2- s;A} )}}{{( p - 1 )!}}=\frac{a_{0}}{{( {s-1/2})}^p}+(-1)^p\sj \binom{j+p-2}{p-1} D\a(j+p-1)(s-1/2)^{j-1}.
\end{align}

Finally, we give a residue theorem which was given by Flajolet and Salvy.
\begin{lem}(\cite{FS1998})\label{lem3.1}
Let $\xi \left( s \right)$ be a kernel function and let $r(s)$ be a rational function which is $O(s^{-2})$ at infinity. Then
\begin{align}\label{2.9}
\sum\limits_{\alpha  \in O} {{\mathop{\rm Res}}{{\left[ {r\left( s \right)\xi \left( s \right)},s = \alpha  \right]}}}  + \sum\limits_{\beta  \in S}  {{\mathop{\rm Res}}{{\left[ {r\left( s \right)\xi \left( s \right)},s = \beta  \right]}}}  = 0.
\end{align}
where $S$ is the set of poles of $r(s)$ and $O$ is the set of poles of $\xi \left( s \right)$ that are not poles $r(s)$ . Here ${\mathop{\rm Re}\nolimits} s{\left[ {r\left( s \right)},s = \alpha \right]} $ denotes the residue of $r(s)$ at $s= \alpha$. The kernel function $\xi \left( s \right)$ is defined by the two requirements: 1. $\xi \left( s \right)$ is meromorphic in the whole complex plane. 2. $\xi \left( s \right)$ satisfies $\xi \left( s \right)=o(s)$ over an infinite collection of circles $\left| s \right| = {\rho _k}$ with ${\rho _k} \to \infty $.
\end{lem}

\section{Evaluations of Euler $T$-sums and Euler $\S$-sums}
Let $B:=\{b_k\}$, $-\infty < k < \infty$ be a sequence of complex numbers with ${b_k} = o\left( {{k^\beta }} \right)\ (\beta  < 1)$ if $k\rightarrow \pm \infty$.
Flajolet and Salvy \cite{FS1998} applied the kernel function
\[\frac{1}{2}\pi \cot( {\pi s} )\frac{{{\psi ^{( {p - 1} )}}( - s )}}{{( {p - 1})!}}\]
to the base function $r(s)=s^{-q}$ to prove every linear sum $S_{p,q}$ whose weight $p+q$ is odd is expressible as a polynomial in zeta values. Next, we replace $\cot(\pi s)\psi^{(p-1)}(-s)$ by $\cot(\pi s;A)\Psi^{(p-1)}(-s;B)$, and use contour integration to evaluate linear (alternating) Euler $T$-sums and Euler $\S$-sums.

\begin{thm}\label{thm1} Let $p>0$ and $q>1$ be positive integers. We have
\begin{align}
&(-1)^{p+q} \su \frac{{\bar N}\B_n(p)}{(n-1/2)^q} a_{n-1}+\su \frac{N\B_n(p)}{(n-1/2)^q}a_n\nonumber\\
&-(-1)^p \sum_{k=0}^{p-1} \binom{p+q-k-2}{q-1} \su \frac{b_n S\a_{n+1}(k+1)}{n^{p+q-k-1}}\nonumber\\
&-b_0\left((-1)^{p+q}\widehat{t}\a(p+q)+\t\a(p+q)\right)\nonumber\\
&+(-1)^p\sum_{k=1}^q \binom{k+p-2}{p-1} D\B(k+p-1)\left((-1)^{q-k}\widehat{t}\a(q-k+1)-\t\a(q-k+1) \right)\nonumber\\&=0.
\end{align}
\end{thm}
\pf We consider the kernel function
\[\pi \cot ( {\pi s};A )\frac{{{\Psi ^{( {p - 1} )}}( 1/2- s;B )}}{{( {p - 1} )!}}\]
and base function $r(s)=(s-1/2)^{-q}$. Clearly, the function $F(s):=\xi(s)r(s)$ only have poles at all integer and $n-1/2$ ($n$ is a positive integer). The only singularities are poles at the integers. At a negative integer $-n$ and positive integer $n$ these two poles are simple and these residues are
\begin{align*}
&{\rm Res}[F(s),s=-n]=(-1)^{p+q} \frac{{\bar N}\B_{n+1}(p)}{(n+1/2)^q}a_n\quad (n\geq 0),\\
&{\rm Res}[F(s),s=n]=\frac{N\B_n(p)}{(n-1/2)^q}a_n \quad (n\geq 1),
\end{align*}
where we used the identities (\ref{b2})-(\ref{2.2}).
From (\ref{2.7}), the pole $n-1/2$ ($n\geq 2$) has order $p$ and the residue is
\begin{align*}
{\rm Res}[F(s),s=n-1/2]=-(-1)^p \sum_{k=0}^{p-1} \binom{p+q-k-2}{q-1} \frac{b_{n-1} S\a_{n}(k+1)}{(n-1)^{p+q-k-1}}.
\end{align*}
From (\ref{2.8}), the pole $1/2$ has order $p+q$ and the residue is
\begin{align*}
&{\rm Res}[F(s),s=1/2]\nonumber\\&=-b_0\left((-1)^{p+q}\widehat{t}\a(p+q)+\t\a(p+q)\right)\nonumber\\
&\quad+(-1)^p\sum_{k=1}^q \binom{k+p-2}{p-1} D\B(k+p-1)\left((-1)^{q-k}\widehat{t}\a(q-k+1)-\t\a(q-k+1) \right).
\end{align*}
Summing these four contributions yields the statement of the theorem. \hfill$\square$

\begin{thm}\label{thm2} Let $p>0$ and $q>1$ be positive integers. We have
\begin{align}
&(-1)^{p+q} \su \frac{{\bar N}\B_{n+1}(p)}{n^q}a_n+\su \frac{N\B_n(p)}{n^q}a_n\nonumber\\
&-(-1)^p\sum_{k=0}^{p-1} \binom{p+q-k-2}{q-1} \su \frac{b_{n-1}S\a_n(k+1)}{(n-1/2)^{p+q-k-1}}\nonumber\\
&+a_0(-1)^p\binom{p+q-1}{q} \widehat{t}\B(p+q)\nonumber\\
&-2(-1)^p\sum_{j=1}^{[q/2]} \binom{p+q-2j-1}{p-1}D\a(2j)\widehat{t}\B(p+q-2j)\nonumber\\
&=0.
\end{align}
\end{thm}
\pf The proof is similar to the previous proof. We consider the kernel function
\[\frac{1}{2}\pi \cot( {\pi s};A)\frac{{{\Psi ^{( {p - 1} )}}( 1/2- s;B )}}{{( {p - 1})!}}\]
and base function $r(s)=s^{-q}$. Then, by a similar argument as in the proof of above, we may easily deduce the desired result.\hfill$\square$

In Theorems \ref{thm1} and \ref{thm2}, setting $A,B\in\{A_1,A_2\}$, by straightforward calculations, we can get the following corollaries.
\begin{cor} For positive integers $p$ and $q>1$,
\begin{align}
&\begin{aligned}\label{4.1}
&(1-(-1)^{p+q}) \su\frac{h^{(p)}_{n-1}}{(n-1/2)^q}\\
=&(-1)^{p+q}\t(p+q)-(-1)^p(1+(-1)^q)\t(p)\t(q)\\
&-(-1)^p\sum_{k=0}^{p-1} ((-1)^k-1)\binom{p+q-k-2}{q-1}\t(k+1)\z(p+q-k-1)\\
&+(-1)^p \sum_{k=1}^{q}(1-(-1)^{q-k}) \binom{k+p-2}{p-1}\t(q-k+1)\z(k+p-1),
\end{aligned}\\
&\begin{aligned}\label{4.2}
&(1+(-1)^{p+q}) \su\frac{{\bar h}^{(p)}_{n-1}}{(n-1/2)^q}\\
=&-(-1)^{p+q}{\bar t}(p+q)+(-1)^p(1+(-1)^q){\bar t}(p)\t(q)\\
&-(-1)^p\sum_{k=0}^{p-1} ((-1)^k+1)\binom{p+q-k-2}{q-1}{\bar t}(k+1)\z(p+q-k-1)\\
&-(-1)^p\sum_{k=1}^q (1+(-1)^{q-k})\binom{k+p-2}{p-1}{\bar t}(q-k+1){\bar \z}(k+p-1),
\end{aligned}\\
&\begin{aligned}\label{4.3}
&(1-(-1)^{p+q}) \su\frac{{\bar h}^{(p)}_{n-1}}{(n-1/2)^q}(-1)^{n-1}\\
=&(-1)^{p+q}\t(p+q)+(-1)^p(1-(-1)^q){\bar t}(p){\bar t}(q)\\
&+(-1)^p\sum_{k=0}^{p-1} ((-1)^k-1)\binom{p+q-k-2}{q-1}\t(k+1){\bar \z}(p+q-k-1)\\
&-(-1)^p \sum_{k=1}^q (1-(-1)^{q-k}) \binom{k+p-2}{p-1}\t(q-k+1){\bar \z}(k+p-1),
&\end{aligned}\\
&\begin{aligned}\label{4.4}
&(1+(-1)^{p+q}) \su\frac{{h}^{(p)}_{n-1}}{(n-1/2)^q}(-1)^{n-1}\\
=&-(-1)^{p+q}{\bar t}(p+q)-(-1)^{p}(1-(-1)^q)\t(p){\bar t}(q)\\
&+(-1)^p\sum\limits_{k=0}^{p-1} ((-1)^k+1)\binom{p+q-k-2}{q-1}{\bar t}(k+1){\bar \z}({p+q-k-1})\\
&+(-1)^p\sum\limits_{k=1}^q  (1+(-1)^{q-k})\binom{k+p-2}{p-1}{\bar t}(q-k+1)\z(p+k-1),
\end{aligned}
\end{align}
where $\z(1):=0$ and $\t(1):=2\log(2)$.
\end{cor}

\begin{cor} For positive integers $p$ and $q>1$,
\begin{align}
&\begin{aligned}\label{4.5}
&(1-(-1)^{p+q})\su \frac{h^{(p)}_n}{n^q}\\
=&-(-1)^p(1+(-1)^q)\t(p)\z(q)-(-1)^p\binom{p+q-1}{p-1}\t(p+q)\\
&-(-1)^p\sum_{k=0}^{p-1} ((-1)^k-1)\binom{p+q-k-2}{q-1}\t(k+1)\t(p+q-k-1)\\
&+2(-1)^p \sum_{j=1}^{[q/2]} \binom{p+q-2j-1}{p-1}\z(2j)\t(p+q-2j)),
\end{aligned}\\
&\begin{aligned}\label{4.6}
&(1+(-1)^{p+q})\su \frac{{\bar h}^{(p)}_n}{n^q}\\
=&(-1)^p(1+(-1)^q){\bar t}(p)\z(q)+(-1)^p\binom{p+q-1}{p-1}{\bar t}(p+q)\\
&-(-1)^p\sum_{k=0}^{p-1} ((-1)^k+1)\binom{p+q-k-2}{q-1}{\bar t}(k+1)\t(p+q-k-1)\\
&+2(-1)^p \sum_{j=1}^{[q/2]} \binom{p+q-2j-1}{p-1}{\bar \z}(2j){\bar t}(p+q-2j)),
\end{aligned}\\
&\begin{aligned}\label{4.7}
&(1+(-1)^{p+q})\su \frac{{\bar h}^{(p)}_n}{n^q}(-1)^{n-1}\\
=&(-1)^p(1+(-1)^q){\bar t}(p){\bar \z}(q)-(-1)^p\binom{p+q-1}{p-1}{\bar t}(p+q)\\
&-(-1)^p\sum_{k=0}^{p-1} ((-1)^k-1)\binom{p+q-k-2}{q-1}{\t}(k+1){\bar t}(p+q-k-1)\\
&+2(-1)^p \sum_{j=1}^{[q/2]} \binom{p+q-2j-1}{p-1}{\z}(2j){\bar t}(p+q-2j)),
\end{aligned}\\
&\begin{aligned}\label{4.8}
&(1-(-1)^{p+q})\su \frac{h^{(p)}_n}{n^q}(-1)^{n-1}\\
=&-(-1)^p(1+(-1)^q)\t(p){\bar \z}(q)+(-1)^p\binom{p+q-1}{p-1}\t(p+q)\\
&-(-1)^p\sum_{k=0}^{p-1} ((-1)^k+1)\binom{p+q-k-2}{q-1}{\bar t}(k+1){\bar t}(p+q-k-1)\\
&+2(-1)^p \sum_{j=1}^{[q/2]} \binom{p+q-2j-1}{p-1}{\bar \z}(2j)\t(p+q-2j),
\end{aligned}
\end{align}
where $\z(1):=0$ and $\t(1):=2\log(2)$.
\end{cor}

\begin{re} Note that when $q=1$, the four identities (\ref{4.3}), (\ref{4.4}), (\ref{4.7}) and (\ref{4.8}) are also hold. We leave the detail to the
interested reader.
\end{re}

Next, we evaluate the quadratic (alternating) Euler $T$-sums and Euler $\S$-sums.

\begin{thm}\label{thm3.5} For positive integers $p,m$ and $q>1$, then
\begin{align}
&(-1)^{p+q+m} \su \frac{{\bar N}\B_n(m){\bar N}\C_n(p)}{(n-1/2)^q} a_{n-1}+\su \frac{N\B_n(m)N\C_n(p)}{(n-1/2)^q}a_n\nonumber\\
&-(-1)^{p+m} \sum_{k=0}^{p+m-1} \binom{p+q+m-k-2}{q-1} \su \frac{b_nc_nS\a_{n+1}(k+1)}{n^{p+q+m-k-1}}\nonumber\\
&-(-1)^m\sum_{j=1}^m \sum_{k=0}^{m-j} \binom{j+p-2}{p-1} \binom{m+q-k-j-1}{q-1} \su \frac{M\C_n(j+p-1)S\a_{n+1}(k+1)b_n}{n^{m+q-j-k}}\nonumber\\
&-(-1)^p\sum_{j=1}^p \sum_{k=0}^{p-j} \binom{j+m-2}{m-1} \binom{p+q-k-j-1}{q-1} \su \frac{M\B_n(j+m-1)S\a_{n+1}(k+1)c_n}{n^{p+q-j-k}}\nonumber\\
&+{\rm Res}[F(s),s=1/2]=0,
\end{align}
where
\begin{align}\label{3.12}
&{\rm Res}[F(s),s=1/2]\nonumber\\&=-b_0c_0\left((-1)^{p+q+m}\widehat{t}\a(p+q+m)+\t\a(p+q+m)\right)\nonumber\\
&\quad+b_0(-1)^p \sum_{j=1}^{m+q} \binom{j+p-2}{j-1} D\C(j+p-1)\nonumber\\&\quad\quad\quad\quad\times\left((-1)^{m+q-j}\widehat{t}\a(m+q-j+1)-\t\a(m+q-j+1)\right)\nonumber\\
&\quad+c_0(-1)^m \sum_{j=1}^{p+q} \binom{j+m-2}{j-1} D\B(j+m-1)\nonumber\\&\quad\quad\quad\quad\times\left((-1)^{p+q-j}\widehat{t}\a(p+q-j+1)-\t\a(p+q-j+1)\right)\nonumber\\
&\quad+(-1)^{p+m} \sum_{j_1+j_2\leq q+1,\atop j_1,j_2\geq 1} \binom{j_1+m-2}{j_1-1}\binom{j_2+p-2}{j_2-1} D\B(j_1+m-1)D\C(j_2+p-1)\nonumber\\
&\quad\quad\quad\quad\times\left((-1)^{q+1-j_1-j_2}\widehat{t}\a(q+2-j_1-j_2)-\t\a(q+2-j_1-j_2)\right).
\end{align}
\end{thm}
\pf We consider the kernel function
\[ \cot( {\pi s};A )\frac{\Psi^{(m-1)}(1/2-s;B){{\Psi ^{( {p - 1} )}}(1/2- s;C )}}{{(m-1)!( {p - 1} )!}}\]
and base function $r(s)=(s-1/2)^{-q}$. It is obvious that the function
\[ F(s):=\cot( {\pi s};A )\frac{\Psi^{(m-1)}(1/2-s;B){{\Psi ^{( {p - 1} )}}(1/2- s;C)}}{{(m-1)!( {p - 1} })!(s-1/2)^{q}}\]
has simple poles at $s=-n\ (n\geq 0)$ with residues
\begin{align*}
{\rm Res}[F(s),s=-n]=(-1)^{p+q+m}\frac{{\bar N}\B_{n+1}(m){\bar N}\C_{n+1}(p)}{(n+1/2)^q} a_{n},
\end{align*}
and simple poles at $s=n\ (n\geq 1)$, with residues
\begin{align*}
{\rm Res}[F(s),s=n]=\frac{N\B_n(m)N\C_n(p)}{(n-1/2)^q}a_n,
\end{align*}
where we used the identities (\ref{2.1}) and (\ref{2.2}). Clearly, $F(s)$ has poles of order $p+m$ at $s=n-1/2\ (n\geq 2)$. Using (\ref{2.4}) and (\ref{2.7}) we find that the residues
\begin{align*}
&{\rm Res}[F(s),s=n-1/2]\\
&=-(-1)^{p+m} \sum_{k=0}^{p+m-1} \binom{p+q+m-k-2}{q-1} \frac{b_{n-1}c_{n-1}S\a_{n}(k+1)}{(n-1)^{p+q+m-k-1}}\\
&\quad-(-1)^m\sum_{j=1}^m \sum_{k=0}^{m-j} \binom{j+p-2}{p-1} \binom{m+q-k-j-1}{q-1}  \frac{M\C_{n-1}(j+p-1)S\a_{n}(k+1)b_{n-1}}{(n-1)^{m+q-j-k}}\\
&\quad-(-1)^p\sum_{j=1}^p \sum_{k=0}^{p-j} \binom{j+m-2}{m-1} \binom{p+q-k-j-1}{q-1} \frac{M\B_{n-1}(j+m-1)S\a_{n}(k+1)c_{n-1}}{(n-1)^{p+q-j-k}}.
\end{align*}
Moreover, $F(s)$ also has a pole of order $p+q+m$ at $s=1/2$. Using (\ref{2.8}) we deduce the formula (\ref{3.12}) by a direct calculation.
Hence, combining these four residue results, we can obtain the desired evaluation. \hfill$\square$

\begin{thm} For positive integers $p,m$ and $q>1$, then
\begin{align}\label{3.13}
&(-1)^{p+q+m} \su \frac{{\bar N}\B_{n+1}(m){\bar N}\C_{n+1}(p)}{n^q} a_{n}+\su \frac{N\B_n(m)N\C_n(p)}{n^q}a_n\nonumber\\
&-(-1)^{p+m} \sum_{k=0}^{p+m-1} \binom{p+q+m-k-2}{q-1} \su \frac{b_{n-1}c_{n-1}S\a_{n}(k+1)}{(n-1/2)^{p+q+m-k-1}}\nonumber\\
&-(-1)^m\sum_{j=1}^m \sum_{k=0}^{m-j} \binom{j+p-2}{p-1} \binom{m+q-k-j-1}{q-1} \su \frac{M\C_{n-1}(j+p-1)S\a_{n}(k+1)b_{n-1}}{(n-1/2)^{m+q-j-k}}\nonumber\\
&-(-1)^p\sum_{j=1}^p \sum_{k=0}^{p-j} \binom{j+m-2}{m-1} \binom{p+q-k-j-1}{q-1} \su \frac{M\B_{n-1}(j+m-1)S\a_{n}(k+1)c_{n-1}}{(n-1/2)^{p+q-j-k}}\nonumber\\
&+{\rm Res}[G(s),s=0]=0,
\end{align}
where
\begin{align}\label{3.14}
&{\rm Res}[G(s),s=0]\nonumber\\&=a_0(-1)^{p+m} \sum_{k_1+k_2=q,\atop k_1,k_2\geq 0} \binom{m+k_1-1}{k_1}\binom{p+k_2-1}{k_2} \widehat{t}\B(m+k_1)\widehat{t}\C(p+k_2)\nonumber\\
&\quad-2(-1)^{m+p}\sum_{j=1}^{[q/2]} \sum_{k_1+k_2=q-2j,\atop k_1,k_2\geq 0} \binom{m+k_1-1}{k_1}\binom{p+k_2-1}{k_2} D\a(2j)\widehat{t}\B(m+k_1)\widehat{t}\C(p+k_2).
\end{align}
\end{thm}
\pf We consider the kernel function
\[\pi \cot( {\pi s};A)\frac{\Psi^{(m-1)}(1/2-s;B){{\Psi ^{( {p - 1} )}}(1/2- s;C )}}{{(m-1)!( {p - 1} )!}}\]
and base function $r(s)=s^{-q}$. By the same calculation as in the proof of Theorem \ref{thm3.5}, we thus immediately deduce (\ref{3.13}) and (\ref{3.14}) to complete the proof. \hfill$\square$

\begin{thm} For positive integers $p,m$ and $q>1$, then
\begin{align}
&(-1)^{p+q+m} \su \frac{{\bar M}\B_n(m){\bar N}\C_{n+1}(p)}{n^q}a_n +\su \frac{M\B_n(m)N\C_n(p)}{n^q}a_n\nonumber\\
&+(-1)^m\sum_{k=0}^m \binom{m+q-k-1}{q-1}\binom{p+k-1}{p-1} \su \frac{N\C_n(p+k)}{n^{m+q-k}}a_nb_n\nonumber\\
&-(-1)^m \sum_{k+j\leq m+1,\atop k,j\geq 1} \binom{m+q-k-j}{q-1}\binom{p+k-2}{p-1}\su \frac{R\a_n(j)N\C_n(p+k-1)}{n^{m+q+1-k-j}}b_n\nonumber\\
&-(-1)^p\sum_{k_1+k_2+k_3=p-1,\atop k_1,k_2,k_3\geq 0} \binom{m+k_2-1}{m-1}\binom{q+k_3-1}{q-1} \su \frac{S\a_n(k_1+1)N\B_n(m+k_2)}{(n-1/2)^{k_3+q}}c_{n-1}\nonumber\\
&+{\rm Res}[H(s),s=0]=0,
\end{align}
where
\begin{align}
&{\rm Res}[H(s),s=0]\nonumber\\
&=a_0b_0(-1)^p \binom{p+q+m-1}{p-1} \widehat{t}\C(p+q+m)\nonumber\\
&\quad+a_0(-1)^{p+m} \sum_{j=1}^{q+1} \binom{j+m-2}{m-1}\binom{p+q-j}{p-1}D\B(j+m-1)\widehat{t}\C(p+q+1-j)\nonumber\\
&\quad-2b_0(-1)^p \sum_{j=1}^{[(m+q)/2]}\binom{p+q+m-2j-1}{p-1} D\a(2j)\widehat{t}\C(p+q+m-2j)\nonumber\\
&\quad-2(-1)^{p+m}\sum_{2j_1+j_2\leq q+1,\atop j_1,j_2\geq 1} \binom{j_2+m-2}{m-1} \binom{p+q-2j_1-j_2}{p-1}\nonumber\\&\quad\quad\quad\quad\quad\quad\quad\quad\times D\a(2j_1)D\B(j_2+m-1)\t\C(p+q+1-2j_1-j_2).
\end{align}
\end{thm}
\pf We consider the kernel function
\[\pi \cot ( {\pi s};A )\frac{\Psi^{(m-1)}(-s;B){{\Psi ^{( {p - 1} )}}(1/2- s;C )}}{{(m-1)!( {p - 1} )!}}\]
and base function $r(s)=s^{-q}$. By direct residue computations, we can obtain the desired evaluation with the help of formulas (\ref{b1})-(\ref{2.4}). \hfill$\square$

It is clear that the main results in our previous paper \cite{XW2019} are immediate corollaries of this paper. Moreover, it is possible that of some other relations involving alternating Euler $T$-sums and related sums can be proved by using the techniques of the present paper. For example, let $A^{(l)}:=\{a_k^{(l)}\}\ (-\infty < k < \infty,\ l$ is any positive integer) be any sequences of complex numbers with ${a_k^{(l)}} = o\left( {{k^\alpha }} \right)\ (\alpha  < 1)$ if $k\rightarrow \pm \infty$,
 consider these two function
\[\cot \left( {\pi s};A \right)\frac{\Psi^{(p_1-1)}(1/2-s;A^{(1)}){{\Psi ^{( {p_2- 1} )}}(1/2- s;A^{(2)})\cdots \Psi^{(p_r-1)}(1/2-s;A^{(r)})}}{{(m-1)!( {p - 1} !})(s-1/2)^{q}}\]
and
\[\cot \left( {\pi s};A \right)\frac{\Psi^{(p_1-1)}(1/2-s;A^{(1)}){{\Psi ^{( {p_2 - 1} )}}(1/2- s;A^{(2)})\cdots \Psi^{(p_r-1)}(1/2-s;A^{(r)})}}{{(m-1)!( {p - 1} !})s^{q}}\]
we can deduce the following results
\begin{align*}
&(-1)^{p_1+\cdots+p_r+q}\su \frac{{\bar N}^{(A^{(1)})}_n(p_1){\bar N}^{(A^{(2)})}_n(p_2)\cdots{\bar N}^{(A^{(r)})}_n(p_r)}{(n-1/2)^q}a_{n-1}\\
&\quad+\su \frac{{ N}^{(A^{(1)})}_n(p_1){ N}^{(A^{(2)})}_n(p_2)\cdots{ N}^{(A^{(r)})}_n(p_r)}{(n-1/2)^q}a_{n}\\
&\quad+\sum (\text{sums of degree}\leq r-1)=0
\end{align*}
and
\begin{align*}
&(-1)^{p_1+\cdots+p_r+q}\su \frac{{\bar N}^{(A^{(1)})}_{n+1}(p_1){\bar N}^{(A^{(2)})}_{n+1}(p_2)\cdots{\bar N}^{(A^{(r)})}_{n+1}(p_r)}{n^q}a_{n}\\
&\quad+\su \frac{{ N}^{(A^{(1)})}_n(p_1){ N}^{(A^{(2)})}_n(p_2)\cdots{ N}^{(A^{(r)})}_n(p_r)}{n^q}a_{n}\\
&\quad+\sum (\text{sums of degree}\leq r-1)=0,
\end{align*}
but we can't get the general explicit formulas.
\\[5mm]
{\bf Acknowledgments.}  The authors express their deep gratitude to Professors Masanobu Kaneko and Jianqiang Zhao for valuable discussions and comments.

 \end{document}